\documentclass[12pt,a4paper]{amsart}

\usepackage{amsmath, amsfonts, xifthen, latexsym, amssymb, amsthm, amscd}

\usepackage[utf8]{inputenc}
\usepackage{graphicx}
\usepackage{url}

\usepackage{hyperref}
\hypersetup{colorlinks=true,citecolor=blue,filecolor=blue,linkcolor=blue,urlcolor=blue}
\usepackage[margin=1.4in]{geometry}

\newtheorem{theorem}{Theorem}
\newtheorem{lemma}{Lemma}[section]

\newtheorem{proposition}{Proposition}[section]

\newtheorem{corollary}{Corollary}[section]

\theoremstyle{definition}
\newtheorem{rules}{Rule}
\newtheorem{definition}{Definition}[section]

\newcommand{\eps}{{\varepsilon}}

\renewcommand{\phi}{{\varphi}}

\newcommand{\Free}{\operatorname{Free}}

\renewcommand{\le}{\leqslant}
\newcommand{\cC}{\mathcal{C}}\newcommand{\cF}{\mathcal{F}}\newcommand{\cH}{\mathcal{H}}

\newcommand\N{\mathbb N}
\newcommand\cG{\mathcal{G}}

\newcommand{\actson}{\curvearrowright}

\newcommand{\cal}[1]{{\mathcal #1}}

\usepackage{etoolbox}
\newtoggle{no_cases}\toggletrue{no_cases}
\newcommand{\case}[2][]{\iftoggle{no_cases}{\left\{\begin{array}{ll}#2 & #1}{\\#2 & #1}\togglefalse{no_cases}}
\newcommand{\esac}{\end{array}\right.\toggletrue{no_cases}}
\newcommand{\hatt}{\hat{\Theta}}
\newcommand{\ovet}{\overline{\Theta}}
\newcommand{\Stab}{\operatorname{Stab}}\newcommand{\Sub}{\operatorname{Sub}}

\begin{document}
\title[Free minimal actions]{Free minimal actions of countable groups with invariant probability measures}
\author{Gábor Elek}
\address{Department of Mathematics And Statistics, Fylde College, Lancaster University, Lancaster, LA1 4YF, United Kingdom}

\email{g.elek@lancaster.ac.uk}  

\thanks{The author was partially supported
by the ERC Consolidator Grant "Asymptotic invariants of discrete groups,
sparse graphs and locally symmetric spaces" No. 648017. }

\begin{abstract} We prove that for any countable group $\Gamma$ there exists a
  free minimal continuous action $\alpha:\Gamma\actson \cC$ on the Cantor set admitting
  an invariant Borel probability measure. 

\end{abstract}\maketitle

\setcounter{tocdepth}{2}
\section{Introduction}
In \cite{HM} Hjorth and Molberg proved that for any countable group $\Gamma$
there exists a free and continuous action of $\Gamma$ on the Cantor set $\cC$ that admits
an invariant Borel probability measure. Our main goal is to prove that one can also
assume the minimality of the action. 
\begin{theorem} \label{fotetel}
For any countably infinite group $\Gamma$, there exists a free minimal continuous action $\alpha: \Gamma \actson \cC$ 
on the Cantor set admitting an invariant Borel probability measure.
\end{theorem}
\noindent
In fact, we show that $\alpha$ can be chosen to be universal in the following sense.
\begin{theorem} \label{fotetel2} There exists a free, minimal action $\alpha:\Gamma\actson\cC$ that satisfies the following condition.
Let $\beta:\Gamma\actson X$ be a free Borel action of the above group $\Gamma$ on
the standard Borel space $X$ admitting an invariant Borel probability measure
$\mu$. Then there exists a Borel embedding $\Phi:X'\to \cC$ commuting with the actions $\alpha$ and $\beta$,  such
that $X'$ is an invariant Borel set of $X$ and  $\mu(X')=1$.
\end{theorem}
\noindent
It is important to note that Weiss proved \cite{WMin} the existence of a 
topologically
free minimal Cantor action satisfying the condition of Theorem
\ref{fotetel2}. Notice however, that topologically free minimal Cantor actions of the free group
admitting an invariant measure can be very far from being essentially free. In fact, it is possible
that the  equivalence relation associated to such topologically free minimal action is hyperfinite \cite{AbE}. 
Our approach to prove Theorem \ref{fotetel2} is based on the philosophy behind Weiss' proof and the notion of
properness. 

\noindent
First we prove a result on Borel universality. 
\begin{proposition} \label{free}
For any countably infinite group $\Gamma$ and for any free Borel action $\beta:\Gamma\actson X$ on the standard Borel space, there exists
an injective equivariant Borel map $\Psi'_\beta: X \to \Free(\cC^\Gamma)$ ( where $\Free(\cC^\Gamma)$ is the free part of the
Bernoulli $\cal{C}$-shift space), such that the closure of the set $\Psi'_\beta(X)$
is still in  $\Free(\cC^\Gamma)$.
\end{proposition}
\noindent
As a consequence, we will obtain the following theorem.
\begin{theorem} \label{fotetel3} 
For any countably infinite group $\Gamma$, there exists a free continuous action $\zeta:\Gamma\actson \cC$ such that
for any free Borel action $\beta:\Gamma\actson X$ on the standard Borel space we have an injective Borel map $\Psi_\beta:X\to \cC$
satisfying $\Psi_\beta\circ \beta=\zeta\circ\Psi_\beta.$
\end{theorem}
\noindent
Note that Seward and Tucker-Drob \cite{ST} (see also \cite{Bern}) 
proved the following result:
For any free Borel action $\alpha:\Gamma\actson X$ there exists a (not necessarily
injective)  equivariant Borel map $\Psi: X\to \Free(\{0,1\}^\Gamma)$
such that the closure of the set $\Psi(X)$ is still in $\Free(\{0,1\}^\Gamma)$. 
\vskip 0.1in
\noindent
In 1952 Oxtoby \cite{Oxt} proved that there exists a free, minimal $Z$-action on the Cantor set that is
not uniquely ergodic, that is, it admits more than one ergodic invariant Borel probability measures.
We will prove the following corollary of Theorem \ref{fotetel2}.
\begin{corollary}\label{teteloxt}
Any countably infinite group
$\Gamma$ has a free, minimal action on the Cantor set that is not uniquely ergodic.
\end{corollary}
\noindent
We will also prove a version of Theorem \ref{fotetel3} for uniformly recurrent subgroups (Theorem \ref{urs}) answering a question of
Glasner and Weiss in the ``universal sense''. 
Finally, we prove Theorem \ref{fotetel}. The main idea is to construct an explicit free minimal action for any countably infinite group using an inductive
``learning'' algorithm. Then
we combine this construction with Theorem \ref{fotetel3} to obtain our main result.

\vskip 0.1in
\noindent
{\bf Acknowledgement:} We are grateful to Benjamin Weiss for sending us his paper \cite{WMin}. 

\section{The proof of Theorem \ref{fotetel3}}\label{ketto}
Let $\Gamma$ be a countable infinite group and $\{\sigma_i\}^\infty_{i=1}$ be a generating system of $\Gamma$. Also for $n\geq 1$, let $\Gamma_n$ be the
subgroup of $\Gamma$ generated by the elements $\{\sigma_i\}^n_{i=1}$. We will also assume that if $\Gamma_n\neq \Gamma$, then $\sigma_{n+1}\notin \Gamma_n$. Hence, 
if $\Gamma_n\neq \Gamma$ then $\Gamma_{n+1}\neq \Gamma_{n}$. Consequently, if $\Gamma_n$ is finite, then $|\Gamma_n|\geq 2^n$.
Let $\alpha:\Gamma\actson X$ be a Borel action of $\Gamma$ on the standard Borel space $X$.
We define a sequence $\{\cG_n\}^\infty_{n=1}$ of Borel graph structures on $X$ in the following way.
If $p,q\in X$, $p\neq q$, then let $(p,q)\in E(\cG_n)$ if there exists $1\leq i \leq n$ such that $\alpha(\sigma_i)(p)=q$ or $\alpha(\sigma_i)(q)=p$.
A Borel $\cal{C}$-coloring of $X$ is a Borel map $\phi:X\to \cC$, where $\cC=\{0,1\}^\N$ is the Cantor set. We say that $\phi$ is a {\bf proper} $\cC$-coloring
with respect to $\alpha:\Gamma\actson X$ if for any $r>0$ there exists $S_r>0$ such that
for any $p,q\in X$ $(\phi(p))_{S_r}\neq (\phi(q))_{S_r}$, provided that $0<d_{\cG_r}(p,q)\leq r$, where
\begin{itemize}
\item $d_{\cG_r}$ is the shortest path metric on the components of the Borel graph $\cG_r$.
\item For $\kappa\in \cC$ and $s>0$, $(\kappa)_s\in \{0,1\}^{[s]}$, denotes the projection of $\kappa$ onto its first $s$ coordinates.
Here $[s]$ denotes the set $\{1,2,\dots,s\}$.
\end{itemize}
\noindent
The Borel coloring $\phi$ is called separating if for any $p\neq q\in X$, $\phi(p)\neq \phi(q)$.
\begin{lemma} \label{proper} For any Borel action $\alpha:\Gamma\actson X$ there exists a separating, proper $\cC$-coloring with respect to $\alpha$.
\end{lemma}
\proof
First, for any $r\geq 1$ we construct
a new Borel graph $\cH_r$ of bounded vertex degree on $X$ such that $(p,q)\in E(\cH_r)$ if $0<d_{G_r}(p,q)\leq r$.
By the classical result of Kechris, Solecki and Todorcevic \cite{KST}, there exists an integer $m_r>0$ and a Borel coloring
$\psi_r:X\to \{0,1\}^{[m_r]}$ such that
$\psi_r(p)\neq \psi_r(q)$, whenever $p$ and $q$ are adjacent vertices in the Borel graph $\cH_r$.
Then $\phi_1(p)=\{\psi_1(p)\psi_2(p)\dots\}\in\{0,1\}^\N$, defines
a proper $\cC$-coloring of $X$ with respect to $\alpha$. Now we use the usual trick to obtain a separating coloring.
Let $\phi_2:X\to \cC$ be a Borel embedding. Then if $\phi_1(p)=\{a_1 a_2 a_3\dots \}$ and
$\phi_2(p)=\{b_1 b_2 b_3\dots \}$ let
$\phi(p)=\{a_1b_1a_2b_2\dots\}$. Clearly, $\phi$ is a separating, proper $\cC$-coloring with respect to $\alpha$. \qed
\vskip 0.1in
\noindent
Now we prove Proposition \ref{free}. Let $\alpha:\Gamma \actson X$ be a free Borel action and let $\phi:X\to \cC$ be a separating, proper $\cC$-coloring.
Consider the Bernoulli shift $\cC^\Gamma$ with the natural
left action
$$L_\delta(\rho)(\gamma)=\rho(\gamma\delta)\,.$$
\noindent
The map $\Psi'_\alpha:X\to \cC^\Gamma$ is defined as usual by
$$\Psi'_\alpha(x)(\gamma)=\phi(\alpha(\gamma)(x))\,.$$
\noindent
Clearly, $\Psi'_\alpha$ is Borel and $\Gamma$-equivariant and since
$\phi$ is separating, $\Psi'_\alpha$ is injective as well. Also, $\Psi'_\alpha(X)\subset \Free(\cC^\Gamma)$,
where $ \Free(\cC^\Gamma)$ is the free part of the Bernoulli shift, that is, the set of elements $\rho\in \cC^\Gamma$ such
that if $L_\delta(\rho)=\rho$, then $\delta=e_\Gamma$.
We need to show that if
$\rho\in\cC^\Gamma$ is in the closure of $\Psi'_\alpha(X)$, then
$L_\gamma(\rho)\neq \rho$, whenever $\gamma\neq e_\Gamma$.
It is enough to see that if $\gamma\neq e_\Gamma$, then
$\rho(\gamma)\neq \rho(e_\Gamma)$.
Let $\lim_{n\to\infty} \Psi'_\alpha(x_n)=\rho\in\cC^\Gamma$.
Also, let $r>0$ such that
\begin{itemize}
\item $\gamma\in\Gamma_r.$
\item $d_{G_r}(e_\Gamma,\gamma)\leq r$, where $d_{G_r}$ is the shortest path metric on the left Cayley graph $\mbox{Cay}(\Gamma_r,\{\sigma_i\}^r_{i=1})$.
\end{itemize}
\noindent
Thus, for any $n\geq 1$,
$$d_{G_r}(\alpha(\gamma) (x_n), x_n)\leq r\,.$$
\noindent
Hence, for any $n\geq 1$,
$$\left(\phi(\alpha(\gamma)(x_n))\right)_{S_r}\neq (\phi(x_n))_{S_r}\,.$$
\noindent
Since $\rho(\gamma)=\lim_{n\to\infty} \phi(\alpha(\gamma)(x_n))$ and
$\rho(e_\Gamma)=\lim_{n\to\infty} \phi(x_n)$, we have that
$\rho(\gamma)\neq \rho(e_\Gamma)\,.$ Hence our proposition follows. \qed
\vskip 0.1in
\noindent
Now we prove Theorem \ref{fotetel3}.
Let $\iota:\Gamma\actson \Free(\cC^\Gamma)$ be
the natural free action and $\Psi'_\iota:\Free(\cC^\Gamma)\to \Free(\cC^\Gamma)$ be the
injective $\Gamma$-equivariant Borel map defined in Proposition \ref{free}.
Let $W$ be the closure of $\Psi'_\iota(\Free(\cC^\Gamma))$ in $\cC^\Gamma$. Thus, $W\subset \Free(\cC^\Gamma)$ is a closed invariant subset.
Then, $W$ can be written as the disjoint union of $W_1$ and $W_2$, where $W_1$ is a countable invariant subset and $W_2$ is a closed invariant subset
homeomorphic to the Cantor set. Indeed, $W$ is an uncountable compact set. Hence $W$ can be written as the disjoint union of
$W_a$ and $W_b$, where $W_a$ is homeomorphic to the Cantor set and $W_b$ is countable. Clearly, the orbit closure $W_2$ of
$W_a$ has no isolated point and $W_1=W\backslash W_2$ is a countable invariant set.
Now, let $\alpha:\Gamma\actson X$ be a free Borel action and $\Psi_\alpha':X\to \Free(\cC^\Gamma)$ be the 
injective $\Gamma$-equivariant  Borel map defined in  Proposition \ref{free}.
Let $\tilde{\Psi}_\alpha=\Psi'_\iota \circ \Psi'_\alpha.$
Clearly, $\tilde{\Psi}_\alpha:X\to W$ is an injective $\Gamma$-equivariant Borel map. Now we define the injective $\Gamma$-equivariant Borel map
$\Psi_\alpha:X\to W_2$ by modifying $\tilde{\Psi}_\alpha$ on countably many orbits. Let $Y=(\tilde{\Psi}_\alpha)^{-1}(W_1).$
We can assume that the complement $W_3$ of $\tilde{\Psi}_\alpha(X)$ in $W_2$ is infinite, so we can define
an injective $\Gamma$-equivariant map $\Psi':Y\to W_3$. 
Therefore, we can define the injective $\Gamma$-equivariant Borel map
$\Psi_\alpha:X\to W_2$ by 
\begin{itemize}
\item $\Psi_\alpha(x)=\tilde{\Psi}_\alpha(x)$ if $x\in X\backslash Y$.
\item $\Psi_\alpha(x)=\Psi'(x)$ if $x\in Y$.
\end{itemize}
\noindent
Since $W_2$ is homeomorphic to $\cC$, our theorem follows. \qed

\section{Uniformly recurrent subgroups} \label{uniform}
Let $\Gamma$ be a countable group and let $\Sub(\Gamma)$ be the  space of subgroups of $\Gamma$ \cite{GW}. Then $\Sub(\Gamma)$ is
a compact, metrizable space and  conjugations define a continuous action $c:\Gamma\actson \Sub(\Gamma)$. Now let $\beta:\Gamma \actson M$ be a continuous minimal nonfree action of a countable group 
on a compact metric space. Then we have a natural equivariant 
map $\Stab_\beta: M\to \Sub(\Gamma)$
from our space $M$ to the compact space of subgroups of $\Gamma$, mapping each point $x\in M$ to its stabilizer subgroup.
Glasner and Weiss (Proposition 1.2 \cite{GW}) proved that the set $M_0$ of the 
continuity points
of the map $\Stab_\beta$ is a dense, invariant $G_\delta$ subset of $M$ and the closure
of $\Stab_\beta(M_0)$ in $\Sub(\Gamma)$ is a minimal closed invariant subset of
$\Sub(\Gamma)$, that is, a {\bf uniformly recurrent subgroup}. They asked if
for any uniformly recurrent subgroup $Z\subset \Sub(\Gamma)$ there exists a
minimal continuous action $\beta:\Gamma \actson  M$ such that $M_0=M$ and $\Stab_\beta(M)=Z$. This question has been
answered in \cite{Elek} and \cite{MT}.
Now we show that one can answer the question of Glasner and Weiss in the
universal sense. 
\begin{theorem} \label{urs} Let $\Gamma$ be a countable group and $Z\subset\Sub(\Gamma)$ be a
uniformly recurrent subgroup. Then there exists a continuous action $\zeta_Z:\Gamma\actson \cC$ such
that
\begin{itemize}
\item $\Stab_{\zeta_Z}: \cC\to \Sub(\Gamma)$ is continuous everywhere and $\Stab_{\zeta_Z}(\cC)=Z$.
\item For any Borel action $\alpha: \Gamma \actson X$  such that for any $x\in X$
the group $\Stab_\alpha(x)$ is in $Z$ (we call these actions $(\Gamma,Z)$-actions),
there exists an injective Borel map $\Psi_\alpha: X \to \cC$ such that $\Psi_\alpha\circ \alpha=\zeta_Z\circ \Psi_\alpha.$
\end{itemize}
\end{theorem}
\noindent
It was proved in Section 5 \cite{Elek} that there exist countable groups $\Gamma$ and uniformly recurrent subgroups
$Z\subset \Sub(\Gamma)$  such that no Borel $(\Gamma,Z)$-action admits an invariant Borel probability measure. However,
we have the following nonfree analogue of the aforementioned result of Hjorth and Molberg.
\begin{corollary}
Let $\Gamma$ be a countable group and let $Z\subset \Sub(\Gamma)$ be a uniformly recurrent subgroup.
If there exists a Borel $(\Gamma,Z)$-action $\alpha:\Gamma\actson X$ that admits an invariant Borel probability measure, then
there exists a continuous $(\Gamma,Z)$-action $\beta:\Gamma\actson \cC$ on the Cantor set admitting an invariant
Borel probability measure such that
\begin{itemize}
\item The map $\Stab_\beta:\cC\to \Sub(\Gamma)$ is continuous everywhere and
\item $\Stab_\beta(\cC)=Z$.
\end{itemize}
\end{corollary}
\proof (of Theorem \ref{urs})
Let $Z\subset \Sub(\Gamma)$ be
a uniformly recurrent subgroup.
We define the Bernoulli shift space $\cC^Z$ of $Z$ in the following way.
Let $$\cC^Z=\cup_{H\in Z} \cal{F}(H)\,,$$
where $\cal{F}(H)$ is the set of maps $\rho:\Gamma/H\to \cC$ from the left coset space of $H$ to the Cantor set.
\noindent
The action of $\Gamma$ on $\cC^Z$ is defined as follows.
\begin{itemize}
\item If $\rho\in \cal{F}(H)$ then
$L_\delta(\rho)\in \cal{F}(\delta H \delta^{-1})$ and
\item $L_\delta(\rho)(\gamma \delta H \delta^{-1})=\rho(\gamma\delta H)\,.$
\end{itemize}
\begin{lemma}
$L:\Gamma\to \mbox{Homeo}(\cC^Z)$ is a homomorphism.
\end{lemma}
\proof
We need to show that if
$\rho:\Gamma/H\to \cC$ and $\delta_1,\delta_2\in \Gamma$, then
$$L_{\delta_1 }(L_{\delta_2}(\rho))=L_{\delta_1\delta_2}(\rho)\,.$$
\noindent
Observe that $$L_{\delta_1} (L_{\delta_2}(\rho))\in\cal{F}(\delta_1\delta_2 H\delta_2^{-1}\delta_1^{-1})\,\,\mbox{and}\,\,
L_{\delta_1\delta_2}(\rho)\in \cal{F}(\delta_1\delta_2H\delta_2^{-1}\delta_1^{-1}).$$
Now
$$L_{\delta_1}(L_{\delta_2}(\rho))(\gamma \delta_1\delta_2H\delta_2^{-1}\delta_1^{-1})= $$ $$
=L_{\delta_2}(\rho)(\gamma\delta_1\delta_2 H \delta_2^{-1})=\rho (\gamma\delta_1\delta_2 H)=
L_{\delta_1\delta_2}(\rho)(\gamma\delta_1\delta_2H\delta_2^{-1}\delta_1^{-1})\,,$$
\noindent
hence our lemma follows. \qed
\vskip 0.1in
\noindent
We can equip $\cC^Z$ with a compact metric structure $d$ such that
$(\cC^Z,d)$ is homeomorphic to the Cantor set and the $\Gamma$-action above is continuous. 
Let $\rho_1:\Gamma/H_1\to \cC$, $\rho_2:\Gamma/H_2\to \cC$ be elements of $\cC^Z$. We say that
$\rho_1$ and $\rho_2$ are $n$-equivalent, $\rho_1 \equiv_n \rho_2$ if
\begin{itemize}
\item For any $\gamma\in \Gamma$, $d_{G_n}(e_\Gamma,\gamma)\leq n$,
$\gamma\in H_1$ if and only if $\gamma\in H_2$. 
\item For any $\gamma\in \Gamma$, $d_{G_n}(e_\Gamma,\gamma)\leq n$,
$$(\rho_1(\gamma H_1))_n=(\rho_2(\gamma H_2))_n\,.$$
\end{itemize}
\noindent
Then we define $d(\rho_1,\rho_2):=\frac{1}{2^n}$ whenever
$\rho_1\equiv_n \rho_2$ and $\rho_1 \not\equiv_{n+1} \rho_2$.
Let $\{H_n\}^\infty_{n=1}, H\in \Sub(\Gamma)$, $\rho\in\cal{F}(H)$ and for
any $n\geq 1$ let $\rho_n\in\cal{F}(H_n)$.
Observe that $\{\rho_n\}^\infty_{n=1}\to \rho$ in the $d$-metric if and only if
\begin{itemize}
\item $H_n\to H$ in the compact space $\Sub(\Gamma)$ and
\item for any $\gamma\in\Gamma$, $\rho_n(\gamma H_n)\to \rho(\gamma H)\,.$
\end{itemize}
\noindent
We can define $\Free(\cC^Z)$ in the usual way.
We have that $\rho:\Gamma/H\to \cC\in \Free(\cC^Z)$ if
$L_\delta(\rho)\neq \rho$ for any $\delta\notin H$.
Clearly, if for any $\delta\notin H$, $\rho(H)\neq \rho(\delta H)$, then $\rho\in \Free(\cC^Z)$. Also, $\Stab_L(\rho)=H$.
Indeed, if $\delta\in H$, then for any $\gamma\in\Gamma$ we have that
$$L_\delta(\rho)(\gamma H)=\rho(\gamma\delta H)=\rho(\gamma H)\,.$$
\noindent 
On the other hand, if $\delta\notin H$, but $\delta H \delta^{-1}=H$, then
$$L_\delta(\rho)(H)=\rho(\delta H)\neq \rho(H)\,.$$
Now let $\alpha:\Gamma\actson X$ be a Borel $(\Gamma,Z)$-action and $\phi:X\to \cC$ be a separating, proper $\cC$-coloring with respect
to $\alpha$.
We define $\Psi^\phi_\alpha:X\to \cC^Z$ as follows.
\begin{itemize}
\item $\Psi^\phi_\alpha(x)\in \cal{F}(H)$, where $H=\Stab_\alpha(x)$.
\item $\Psi^\phi_\alpha(x)(\gamma H)=\phi(\alpha(\gamma)(x))\,.$
\end{itemize}
\noindent
Clearly, $\Psi^\phi_\alpha:X\to \cC^Z$ is an injective Borel map and $\Psi^\phi_\alpha(X)\subset \Free(\cC^Z)$.
\begin{lemma} 
The map $\Psi_\alpha^\phi:X\to \Free(\cC^Z)$ is $\Gamma$-equivariant.
\end{lemma}
\proof
Let $\delta\in\Gamma$. Then
$$L_\delta(\Psi^\phi_\alpha(x))(\gamma \delta H \delta^{-1} )=\Psi^\phi_\alpha(x)(\gamma\delta H)=
\phi(\alpha(\gamma\delta)(x))\,.$$
\noindent
On the other hand,
$$\Psi^\phi_\alpha(\alpha(\delta)(x))(\gamma\delta H\delta^{-1})=\phi(\alpha(\gamma)\alpha(\delta)(x))=
\phi(\alpha(\gamma\delta)(x))\,.\quad\qed $$
\noindent
Now we prove the nonfree analogue of Proposition \ref{free}.
\begin{proposition} \label{free2}
For any countably infinite group $\Gamma$ and for any free Borel 
$(\Gamma,Z)$-action $\alpha:\Gamma\actson  X$, there exists
an injective equivariant Borel map $\Psi'_\alpha: X \to \Free(\cC^Z)$, such that the closure of the set $\Psi'_\alpha(X)$
is still in  $\Free(\cC^Z)$.
\end{proposition}
\proof Let $\phi$ and $\Psi^\phi_\alpha$ be as above.
Let $\{x_n\}^\infty_{n=1}\subset X$ such that
$$\lim_{n\to \infty}\Psi^\phi_\alpha(x_n)=\rho\in\cal{F}(H)$$ \noindent and $\delta \notin H$. We need to show that $\rho\in \Free(\cC^Z)$. 
Observe that \\ $\{\Stab_\alpha(x_n)\}^\infty_{n=1}\to H$ in $\Sub(\Gamma)$. 
Hence, there exists $N>0$ such that $\delta\notin H_n$ if $n\geq N$. By properness, there exists $m>0$ such
that for all $n\geq N$
$$(\phi(\alpha(\delta) (x_n)))_m\neq (\phi(x_n))_m\,.$$
\noindent
Since, $\lim_{n\to\infty} \phi(\alpha(\delta)( x_n))=\rho(\delta H)$ and $\lim_{n\to\infty} \phi(x_n)=\rho(H)$
we have that $\rho(H)\neq \rho(\delta H)$. Hence, we have that $\rho\in \Free(\cC^Z)$. \qed
\vskip 0.1in
\noindent
Now Theorem \ref{urs} follows from Proposition \ref{free2} exactly
the same way as Theorem \ref{fotetel} follows from Proposition
\ref{free}. Note that if $Y\subset \Free(\cC^Z)$ is a closed, invariant subset, then the continuity of $\Stab_L$ on
$Y$ follows from the definition.  \qed

\section{Minimal Cantor labelings I.} \label{negy}
The goal of the next three sections is to present an explicit contstruction of  free
minimal actions for countably infinite groups. Before getting further let us recall some notions from graph theory.
Let $G$ be a graph, $S\subset V(G)$ be a subset and $r,s>0$ be integers.
Then $S$ is an {\bf $s$-net} if for any $x\in V(G)$ there exists
$y\in S$ such that $d_G(x,y)\le s$. Also,
$S$ is called an {\bf $r$-sparse set} if
for any $x\neq y\in S$, $d_G(x,y)>r$.
\begin{definition}
$T\subset V(G)$ is an $S$-maximal $r$-sparse set if $T$ is maximal among sets that satisfy the following two properties.
\begin{itemize}
\item $T\subseteq S$
\item $T$ is $r$-sparse.
\end{itemize}
\end{definition}
\noindent
Note that an $S$-maximal $r$-sparse set is not necessarily a maximal $r$-sparse set.
However, if $T$ is an $S$-maximal $r$-sparse set, then $T$ is also an $S'$-maximal $r$-sparse
set, provided that $T\subseteq S'\subseteq S$. 
\begin{lemma}\label{net}
If $T$ is an $S$-maximal $r$-sparse set, where $S$ is an $s$-net, then
$T$ is an $s+r$-net.
\end{lemma}
\proof
Let $x\in V(G)$. Then, we have $y\in S$ such that $d_G(x,y)\leq s$. By maximality, if
$y\notin T$ then there exists $z\in T$ such that $d_G(y,z)\leq r$. Hence,
$d_G(x,T)\leq r+s$. \qed
\vskip 0.2in
\noindent
Now let $\Gamma$, $\{\sigma_i\}^\infty_{i=1}$, $\{\Gamma_r\}^\infty_{r=1}$ be as above and for $r\geq 1$
let the graph $G_r$ be defined in the following way.
\begin{itemize}
\item $V(G_r)=\Gamma$.
\item $(p,q)\in E(G_r)$ if for some $1\leq i \leq r$, $\sigma_ip=q$ or $\sigma_iq=p$.
\end{itemize}
 We
construct inductively  sequence of positive integers
$$s_1 <r_1 < s_2 <r_2< \dots\,,$$ $$f(1)<f(2)<f(3)<\dots$$
and 
finite sets $\{F_m \}^\infty_{m=1}$ with distinguished elements $q_m\in F_m$
satisfying certain rules. In the following subsection we describe the rules.
\subsection{The Rules}\label{rules}
We choose the numbers $\{s_i\}^\infty_{i=1}$, $\{r_i\}^\infty_{i=1}$, $\{f(i)\}^\infty_{i=1}$ and \\ $\{|F_m|\}^\infty_{m=1}$ in the following order. 
$$s_1, f(1), r_1, |F_1|, s_2, f(2), r_2, |F_2|, s_3, \dots$$
\noindent
Before stating the rules we need a lemma.
\begin{lemma} \label{kival}
For any $T\geq 1$, there exists $n_T>0$ such that if $n_T\leq n$, then there is an element $\gamma\in\Gamma_n$
so that $d_{G_n}(\gamma,e_\Gamma)\geq T$.
\end{lemma}
\proof
Observe that for all $m\geq 1$, $|B_T(G_m,e_\Gamma)|\leq (3m)^T.$ On the other hand, by our definition $|\Gamma_m|\geq 2^m$, hence
our lemma follows.\qed
\vskip 0.1in
Let $s_1=10$ and then define $f(1)$ in such a way that
$$10 |B_{50}(G_{f(1)},e_\Gamma)|<|\Gamma_{f(1)}|\,.$$
\noindent
Note the if $\Gamma_1$ is infinite we can define $f(1)$ being equal to $1$. Then, pick $r_1$ so that
$$ |B_{\frac{1}{10} r_1}(G_{f(1)}, e_\Gamma)|\geq 10 |B_{50}(G_{f(1)},e_\Gamma)|\,.$$
\noindent
Finally, choose $|F_1|$ in such a way that  $|B_{r_1}(G_{f(1)},e_\Gamma)|< |F_1|\,.$
\begin{rules}\label{rule1}
Suppose that $s_1, f(1), r_1, |F_1|, s_2, f(2),\dots, |F_m|$ have already been chosen.
Let $\kappa_m=|F_m|^{|B_{s_m}(\Gamma_{f(m)},e_\Gamma)|}$ be the number of labelings
of the ball $B_{s_m}(\Gamma_{f(m)},e_\Gamma)$ by the set $F_m$.
Choose the number $s_{m+1}$ so large (see Lemma \ref{largeenough} below) so that for any $n\geq n_{1000 r_m \kappa_m}$ (see Lemma \ref{kival})
and for any $2r_m$-net $T\subset \Gamma_n$ one can pick a subset $L_T$ satisfying
the following three conditions.
\begin{itemize}
\item $|L_T|=\kappa_m$.
\item For any $\delta\in L_T$, 
$$\frac{1}{3}s_{m+1}\leq d_{G_{n}}(e_\Gamma,\delta) \leq \frac{2}{3} s_{m+1}\,.$$
\item For any $\gamma\neq\delta \in L_T$,
$$d_{G_n}(\gamma,\delta)> 20 r_m\,.$$
\end{itemize}
\begin{lemma}\label{largeenough}
The number $s_{m+1}$ can be chosen as it is required above.
\end{lemma}
\proof
Let $s_{m+1}=1000 r_m\kappa_m\,.$
Let $c_1< c_2<\dots< c_{\kappa_m}$ be
integers such that for any $1\leq i \leq \kappa_m$ we have
$\frac{1}{3}s_{m+1}+10r_m\leq c_i \leq \frac{2}{3}s_{m+1}-10r_m$
and $|c_{i+1}-c_i|>100 r_m$.
Now using Lemma \ref{kival}, we pick elements $\{\gamma_i\}_{i=1}^{\kappa_m}\subset \Gamma$ such that
$d_{G_n}(e_\Gamma,\gamma_i)=c_i$. For $1\leq i \leq \kappa_m$, let $\delta_i\in T$ such 
that $d_{G_n}(\gamma_i,\delta_i)\leq 2r_m$. Then $L_T:=\{\delta_i\}^{\kappa_m}_{i=1}$ satisfies
the three conditions in Rule \ref{rule1}. \qed

\end{rules}

\begin{rules}\label{rule2}
Suppose that the numbers \\ $s_1, f(1), r_1, |F_1|, s_2, f(2), r_2, \dots, s_{m+1}$ have already been chosen.\\
Choose $f(m+1)$ in such a way that $f(m+1)>f(m)$, $f(m+1)>n_{1000r_m \kappa_m}$ and 
$$10^{m+1}|B_{5s_{m+1}}(G_{f(m+1)},e_\Gamma)|<|\Gamma_{f(m+1)}|\,.$$
\noindent
Then pick the number $r_{m+1}$ so that both conditions below are satisfied.
\begin{itemize}
\item $r_{m+1}\geq 1000 (\sum^{m+1}_{j=1} 5r_{j-1}+ s_j)$
\item $|B_{\frac{1}{10} r_{m+1}}(G_{f(m+1)}, e_\Gamma)|\geq 10^{m+1} |B_{5s_{m+1}}(G_{f(m+1)},e_\Gamma)|\,.$
\end{itemize}
\end{rules}
\begin{rules}\label{rule3}
Suppose that the numbers $s_1, f(1), r_1, |F_1|,\dots, r_{m+1}$ have already been chosen.
Then pick the number $|F_{m+1}|$ in such a way that
$$|B_{r_{m+1}}(G_{f(m+1)},e_\Gamma)|< |F_{m+1}|\,.$$
\noindent
In particular, one has a function $\phi:\Gamma\to F_m$ such that
if $0<d_{G_{f(m)}}(x,y)\leq r_m$ then $\phi(x) \neq \phi(y)\,.$
\end{rules}
\noindent
Note that if the group $\Gamma$ is torsion-free then for any $n\geq 1$, we can choose $f(n)=n$. On the other hand, if $\Gamma$ is locally finite, then all
the subgroups $\Gamma_n$ are finite and we will need to choose all the values $f(n)$ as above. 
\subsection{Clean labelings}
We assume that $\Gamma$ is equipped with a labeling
$$\Theta=\Theta_{\cC}\times \prod_{m=1}^\infty \Theta_m:\Gamma\to
\cC\times \prod_{m=1}^\infty F_m$$
satisfying the following conditions.
\begin{enumerate}
\item $(\Theta_1)^{-1}(q_1)\subset \Gamma$ is a maximal
$r_1$-sparse subset in the graph $G_1$.
\item For any $m> 1$, $(\Theta_m)^{-1}(q_m)\subset \Gamma$ is a $(\Theta_{m-1})^{-1}(q_{m-1})$-maximal \\
$r_m$-sparse subset in the graph $G_{f(m)}$.
\item For any $m\geq 1$, $\Theta_m(x)\neq \Theta_m(y)$ whenever
  $0<d_{G_{f(m)}}(x,y)\leq r_m$.
\item  For any $m\geq 1$,
$$d_{G_{f(m)}}(e_\Gamma, (\Theta_m)^{-1}(q_m))>10 s_m\,.$$
\end{enumerate}
\noindent
We call such a labeling a {\bf clean labeling} of $\Gamma$. It is easy to see that such clean labelings exist
using Rule \ref{rule3}. In Section \ref{end} we will call a labeling {\bf almost clean labeling} if it satisfies the first
three conditions above. 
From now on we use $\cF$ as a shorthand for $\cC\times \prod_{m=1}^\infty F_m$.
Let $\Theta'=\Theta'_{\cC}\times \prod_{m=1}^\infty\Theta'_m:\Gamma\to\cF$ be
a clean labeling and $j\geq 1$ be an integer.
Let $\cF_j=\{0,1\}^{\{1,2,\dots,j\}}\times \prod^j_{m=1} F_j$, let $z\in\Gamma$ and $t>1$ be
an integer. Then $B_t^{\Theta',j}(G_{f(j)},z)$ is a $\cF_j$-labeled $G_{f(j)}$-ball around $z$
of radius $t$, where the label of $y\in B_t(G_{f(j)},z)$ is given as
$$\Theta'_{[j]}(y)=(\Theta'_\cC(y))_{j}\times \prod^{j}_{m=1}
\Theta'_m(y)\,.$$
\subsection{Patchings}
{\bf Patching} is an elementary construction that turns one \\ clean labeling into
another clean labeling using a third clean labeling. 
First, we define ``regular $n$-patchings''.
Suppose that two clean labelings $\Theta^a=\Theta^a_\cC \times
\prod_{m=1}^\infty\Theta^a_m:\Gamma\to\cF$ and $\Theta^d=\Theta^d_\cC \times
\prod_{m=1}^\infty\Theta^d_m:\Gamma\to\cF$  have already been given.
Let $x\in \Gamma$ such that $\Theta^d_n(x)=q_n.$ Consider
the $\cF_n$-labeled $G_{f(n)}$-ball $B=B^{\Theta^d,n}_{s_n+ r_{n-1}}(G_{f(n)},x)$.
\noindent
The ball is the {\bf patch} we wish to insert into $\Gamma$
in order to turn $\Theta^a$ into the new clean labeling $\Theta^b$. 
In the course of the paper we will refer
to the change of labeling described in the next proposition as ``patching the
$\cF_n$-labeled $G_{f(n)}$-ball $B^{\Theta^d,n}_{s_n+r_{n-1}}(G_{f(n)},x)$ into $\Theta^a$
around the vertex $y$.''

\begin{proposition}\label{patching}
Let $y\in\Gamma$ such that $\Theta^a_n(y)=q_n$.
Then, there exists a clean labeling $\Theta^b=\Theta^b_\cC \times
\prod_{m=1}^\infty\Theta^b_m:\Gamma\to\cF$ such that
\begin{itemize}
\item For any $z\in\Gamma$,$(\Theta^b_\cC)_{n^+}=(\Theta^a_\cC)_{n^+}$, where
  for $c\in\{0,1\}^{\N}$, $(c)_{n^+}$ denotes the projection of $c$ onto
  $\{0,1\}^{\N\backslash \{1,2,\dots,n\}}$.
\item For any $z\in\Gamma$ such that
$d_{G_{f(n)}}(y,z)\geq s_n+4r_{n-1}$, we have
  $\Theta^b_{\cC}(z)=\Theta^a_{\cC}(z)$ (in particular, $\Theta^b_{\cC}(z)=\Theta^a_{\cC}(z)$, if $z$ is not in the $\Gamma_{f(n)}$-orbit\\ of $y$).
\item For any $z\in\Gamma$ such that $d_{G_{f(n)}}(y,z)\leq s_n+r_{n-1}$, we have 
$(\Theta^b_{\cC}(z))_{n}=(\Theta_{\cC}^d(zy^{-1}x)_{n}$.
\item For $m\geq n+1$, we have $\Theta^b_m=\Theta^a_m$.
\item For $m\leq n$ and $z\in\Gamma$ such that
$d_{G_{f(n)}}(y,z)\geq s_n+4r_{n-1}$, we have $\Theta^b_m(z)=\Theta^a_m(z)$ (again, $\Theta^b_{m}(z)=\Theta^a_{m}(z)$, if $z$ is not in the $\Gamma_{f(n)}$-orbit of $y$).
\item For $m\leq n$ and $z\in\Gamma$ such that
$d_{G_{f(n)}}(y,z)\leq s_n+r_{n-1}$, we have 
$\Theta^b_m(z)=\Theta_m^d(zy^{-1}x)$.
\end{itemize}
\end{proposition}
\proof
We need to define only $\Theta_m^b:\Gamma\to F_m$ for each $m\leq n-1$. 
In the first round we will define
$(\Theta^b_m)^{-1}(q_m)$ for each $m$, in the second (much easier) round we
will define the remaining part of the maps $\Theta^b_m$.
For $1\leq i \leq n-1$, let
$$\kappa_i\subset B_{s_n+2 r_{n-1}} (G_{f(n)},y)\cup (\Gamma\backslash 
B_{s_n+3r_{n-1}}(G_{f(n)},y))$$
be defined in the following way.
\begin{itemize}
\item If $d_{G_{f(n)}}(z,y)> s_n+3 r_{n-1}$ then let
$z\in \kappa_i$ if and only if $\Theta_i^a(z)=q_i\,.$
\item If $d_{G_{f(n)}}(z,y)\leq s_n+ 2 r_{n-1}$, then let
$z\in \kappa_i$ if and only if $\Theta^d_i(zy^{-1}x)=q_i\,.$
\end{itemize}
\noindent
We also define $\kappa_0$ as the empty set.
One can observe that
$\kappa_{n-1}\subset \kappa_{n-2}\subset\dots\subset \kappa_1$
and for any $1\leq i \leq n-1$, the set $\kappa_i$ is $r_i$-sparse in $G_{f(i)}$.
 Let $\kappa_1\subset\lambda_1$ be an arbitrary maximal $r_1$-sparse set. Now, let 
$\kappa_2\subset \lambda_2$ be a $\lambda_1$-maximal $r_2$-sparse set. Inductively, we construct
the sets $$\lambda_{n-1}\subset\lambda_{n-2}\subset\dots \subset\lambda_1$$
\noindent 
so that for any $1\leq i \leq n-1$ we have $\kappa_i\subset \lambda_i$ and $\lambda_i$ is a $\lambda_{i-1}$-maximal $r_i$-sparse set.
\begin{lemma}
Let $1\leq i \leq n-1$ and $z_i\in B_{s_n+r_{n-1}}(G_{f(n)},y)$.
Also, let $z_i\in \lambda_i$. Then $z_i\in\kappa_i$. Similarly, if $z_i\in \Gamma\backslash B_{s_n+4r_{n-1}}(G_{f(n)},y)$ and $z_i\in\lambda_i$, then
$z_i\in\kappa_i$. 
\end{lemma}
\proof
Suppose that $z_1\in B_{s_n+r_{n-1}}(G_{f(n)},y)$ such that $z_1\in \lambda_1\backslash \kappa_1$. \\
Since $(\Theta^d_1)^{-1}(q_1)$ is a maximal $r_1$-sparse set, there exists
$w_1\in (\Theta^d_1)^{-1}(q_1)$ such that
$d_{G_{f(1)}}(z_1y^{-1}x, w_1)\leq r_1$.
That is, $d_{G_{f(1)}}(z_1,w_1x^{-1}y)\leq r_1$ and \\ $d_{G_{f(1)}}(y,w_1x^{-1}y)\leq s_n+2r_{n-1}.$
Thus, $w_1x^{-1}y\in\kappa_1\subset \lambda_1$ leading to a contradiction.
Suppose that our statement holds for all $1\leq k\leq i$. Let $z_{i+1}\in B_{s_n+r_{n-1}}$ and $z_{i+1}\in\lambda_{i+1}\backslash \kappa_{i+1}$.
By the induction argument, $z_{i+1}\in \kappa_i$. Since $(\Theta^d_{i+1})^{-1}(q_{i+1})$ is a $(\Theta^d_i)^{-1}(q_i)$-maximal $r_{i+1}$-sparse set,
there exists $w_{i+1}\in (\Theta^d_{i+1})^{-1}(q_{i+1})$ such that $d_{G_{i+1}}(z_{i+1}y^{-1}x,w_{i+1})\leq r_{i+1}$.
That is, \\ $d_{G_{i+1}}(z_{i+1},w_{i+1}x^{-1}y)\leq r_{i+1}$ and $d_{G_{i+1}}(y,w_{i+1}x^{-1}y)\leq s_n+2r_{n-1}.$
Therefore, $w_{i+1}x^{-1}y\in\lambda_{i+1}$ in contradiction with the fact that $\lambda_{i+1}$ is $r_{i+1}$-sparse. The proof of the
second part can be done similarly. \qed
\vskip 0.1in
\noindent
Now, for $1\leq i \leq n-1$ let
\begin{itemize}
\item $\Theta^b_i(z)=q_i$ if $z\in\lambda_i$.
\item $\Theta^b_i(z)=\Theta^d_i(zy^{-1}x)$ if $z\in B_{s_n+r_{n-1}}(G_{f(n)},y).$
\item  $\Theta^b_i(z)=\Theta^a_i(z)$ if $z\in (\Gamma\backslash B_{s_n+4r_{n-1}}(G_{f(n)},y)).$
\end{itemize}
\vskip 0.1in
\noindent
Also, let $\Theta^b_n(z)=\Theta^d_n(zy^{-1}x)$ if $z\in B_{s_n+r_{n-1}}(G_{f(n)},y)$ and
$\Theta^b_n(z)=\Theta^a_n(z)$ if $z\in (\Gamma\backslash B_{s_n+4r_{n-1}}(G_{f(n)},y)).$
Then for $1\leq i \leq n$, let us extend $\Theta^b_i$ onto the set
$$\{z\,\mid s_n+r_{n-1} < d_{G_{f(n)}}(y,z) \leq s_n+4r_{n-1})\}\backslash (\Theta_i^b)^{-1}(q_i)$$
\noindent
to obtain a clean labeling.
Note that $1_\Gamma\notin B_{s_n+ 5r_{n-1}}(G_{f(n)},y)$, hence
by Rule \ref{rule2} such extensions
$\Theta^b_i$ clearly exist. \qed
\vskip 0.2in
We need another kind of patching in our construction that we call ``supersize
$n$-patchings''.
Let $x\in\Gamma$ such that $\Theta^d_n(x)=q_n$. Consider the $\cF_n$-labeled
$G_{f(n)}$-ball $B= B^{\Theta^d,n}_{3r_n}(G_{f(n)},x)$.
The following proposition is about how to insert $B$ into $\Theta^a$.
Notice that in the case of regular patchings the ball $B$ we use as a patch
contained only one element of $(\Theta_n^d)^{-1}(q_n)$.
In the case of supersize patchings the ball $B$ contains many elements of
$(\Theta_n^d)^{-1}(q_n)$. Nevertheless, $B$ will not contain any element of
the set $(\Theta_{n+1}^d)^{-1}(q_{n+1})$. The next proposition is the
``supersized'' version of Proposition \ref{patching}.
\begin{proposition}\label{supersize}
Let $y\in\Gamma$ such that $\Theta^a_n(y)=q_n$ and \\
$d_G(y,(\Theta^a_{n+1})^{-1}(q_{n+1}))> 20 r_n$.
Then, there exists a clean labeling $\Theta^b=\Theta^b_\cC \times
\prod_{m=1}^\infty\Theta^b_m:\Gamma\to\cF$ such that
\begin{itemize}
\item For any $z\in\Gamma$,$(\Theta^b_\cC)_{n^+}=(\Theta^a_\cC)_{n^+}$.
\item For any $z\in\Gamma$ such that
$d_G(y,z)\geq 7r_{n}$, we have
  $\Theta^b_{\cC}(z)=\Theta^a_{\cC}(z)$.
\item For any $z\in\Gamma$ such that $d_G(y,z)\leq 3r_{n}$, we have \\
$(\Theta^b_{\cC}(z))_n=(\Theta_{\cC}^d(zy^{-1}x))_n$.
\item For $m\geq n+1$, let $\Theta^b_m=\Theta^a_m$.
\item For $m\leq n$ and $z\in\Gamma$ such that
$d_{G_{f(n)}}(y,z)\geq 7r_{n}$, we have $\Theta^b_m(z)=\Theta^a_m(z).$
\item For $m\leq n$ and $z\in\Gamma$ such that
$d_{G_{f(n)}}(y,z)\leq 3r_{n}$, we have 
$\Theta^b_m(z)=\Theta_m^d(zy^{-1}x)$.
\end{itemize}
\end{proposition}
\proof
First we construct sets $\lambda_n\subset\lambda_{n-1}\subset\dots\subset\lambda_1$, where $\kappa_i\subset \lambda_i$ for any $1\leq i \leq n$,
$\lambda_1$ is a maximal $r_1$-sparse set and if $i>1$, then
$\lambda_i$ is a $\lambda_{i-1}$-maximal $r_i$-sparse set. Then we can proceed in the same way as in the proof
of Proposition \ref{patching}. \qed
\section{Minimal Cantor labelings II. (codeballs)} \label{codeball}

Let $\Theta'=\Theta'_{\cC}\times \prod_{m=1}^\infty\Theta'_m:\Gamma\to\cF$ and 
$\Theta''=\Theta''_{\cC}\times \prod_{m=1}^\infty\Theta''_m:\Gamma\to\cF$ be clean labelings. Let $\rho>0$ be an integer. 
We say that the $\cF_j$-labeled $G_{f(j)}$-ball $B_t^{\Theta',j}(G_{f(j)},z)$ is {\bf  $(\rho,j)$-syndetic}
in $\Theta''$ if for any $g\in\Gamma$ there exists $h\in\Gamma$ such that
\begin{itemize}
\item $d_{G_{f(j)}}(g,h)\leq \rho$.
\item $B_t^{\Theta',j}(G_{f(j)},z)$ is isomorphic (as $\cF_j$-labeled $G_{f(j)}$-balls) to \\
$B_t^{\Theta'',j}(G_{f(j)},h)$, \\ $B_t^{\Theta',j}(G_{f(j)},z)\cong B_t^{\Theta'',j}(G_{f(j)},h)$,
that is for any $y\in B_t(G_{f(j)},z)$ we have that
$$\Theta'_{[j]}(y)=\Theta''_{[j]}(yz^{-1}h)\,.$$
\end{itemize}
\noindent
\begin{definition}
Let $\Theta'$ and $\Theta''$ be labelings as above. We say that the
$\cF_j$-labeled $G_{f(j)}$-ball $B_l^{\Theta'',j}(G_{f(j)},x)$ {\bf fully contains} an isomorphic copy of
the $G_{f(j)}$-ball $B_{l'}^{\Theta',j}(G_{f(j)},x')$ if there exists $w\in\Gamma$,
$d_{G_{f(j)}}(x,w)+l'<l$ such that $B_{l'}^{\Theta',j}(G_{f(j)},x')$ is isomorphic to
$B_{l'}^{\Theta'',j}(G_{f(j)},w)$ as $\cF_j$-labeled $G_{f(j)}$-balls.
\end{definition}

\noindent
We will start with a clean labeling $\Theta$ and inductively construct a sequence
of clean labelings $\{\Theta^i\}^\infty_{i=1}$ such that
$\Theta^1=\Theta$ and for any $z\in\Gamma$ there exists an integer $i_z>0$
such that
if $i\geq i_z$ then $\Theta^i(z)=\Theta^{i_z}(z)\,.$
Hence, for each $z\in\Gamma$, $\lim_{i\to\infty}\Theta^i(z)=\Theta^\infty(z)$
exists.
The clean labeling $\Theta^\infty:\Gamma\to\cF$ will generate a free, minimal
subshift in the Bernoulli shift space $\cF^\Gamma$.
\vskip 0.1in
\noindent
{\bf The construction of $\Theta^2$}. Our input is the clean labeling $\Theta=\Theta^1$.
Pick any $z_2\in\Gamma$ such that $\Theta_2^1(z_2)=q_2$.
First we modify the labeling $\Theta^1$ on the $G_{f(2)}$-ball
$B_{s_2}(G_{f(2)},z_2)$.
Let $\cal{A}^1=\{A^1_i\}_{i=1}^{\tau_1}$
be the set of all $\cF_1$-labeled $G_{f(1)}$-balls of radius $s_1$ in $\Theta^1$ up to
isomorphism. Let $\{w^1_i\}_{i=1}^{\tau_1}\subset B_{s_2}(G_{f(2)},z_2)\cap (\Theta_1^1)^{-1}(q_1)$ be a subset of vertices
such that
\begin{itemize}
\item For any $1\leq i\leq \tau_1$,
$$\frac{1}{3} s_2\leq d_{G_{f(2)}}(z_2,w^1_i)\leq \frac{2}{3} s_2\,.$$
\item For any $1\leq i\neq j\leq \tau_1$,
$$d_{G_{f(2)}}(w^1_i,w^1_j)>20 r_1\,.$$
\end{itemize}
\noindent
Then for each $1\leq i \leq \tau_1$
pick $z^1_i\in (\Theta_1^1)^{-1}(q_1)$ such that
$B_{2r_1}^{\Theta^1,1}(G_{f(1)},z^1_i)$ fully contains an isomorphic copy of the $G_{f(1)}$-ball $A^1_i$. 
Then using Proposition \ref{supersize}, for each $1\leq i \leq \tau_1$ let us simultaneously patch the
balls $B_{3r_1}^{\Theta^1,1}(G_{f(1)},z^1_i)$ into $\Theta^1$ around $w^1_i$. Since
$(\Theta_1^1)^{-1}(q_1)$ is a maximal $r_1$-sparse set such a family of vertices
$\{z^1_i\}^{\tau_1}_{i=1}$ indeed exists.
Hence, we obtain a modified clean labeling $\hat{\Theta}^2:\Gamma\to\cF$.
We call the $\cF_2$-labeled $G_{f(2)}$-ball $B_2=B_{s_2}^{\hat{\Theta}^2,2}(G_{f(2)},z_2)$ the
{\bf 2-codeball}. We finish the construction of $\Theta^2$ by patching \\
the ball $B_{s_2+r_1}^{\hat{\Theta}^2,2}(G_{f(2)},z_2)$ simultaneously into $\Theta^1$ centered 
around all vertices
$y\in(\Theta_2^1)^{-1}(q_2)$, using Proposition \ref{patching}.
Since $(\Theta^1_1)^{-1}(q_1)$ is an $r_1$-net and 
$(\Theta^1_2)^{-1}(q_2)$ is a $(\Theta^1_1)^{-1}(q_1)$-maximal $r_2$-sparse net,
by Lemma \ref{net} each ball $A^1_j\in\cal{A}^1$ is  
$(2r_2,2)$-syndetic in $\Theta^2$.
Note that the $2$-codeballs are the analogues of the $2$-welcome words of
\cite{WMin}. We will see that for any $i\geq 1$ and 
$x\in (\Theta_2^i)^{-1}(q_2),$
$$B_{s_2}^{\Theta^i,2}(G_{f(2)},x)\cong B_2\,.$$
Therefore, all balls $A^1_j\in\cal{A}^1$ remain
$(2r_2,2)$-syndetic in the process. It is important to note
that it is possible that the set of $\cF_1$-labeled $G_{f(1)}$-balls of radius
$r_1$ will increase when we turn the clean labeling $\Theta^1$
into $\Theta^2$.
\vskip 0.1in
\noindent{\bf The construction of $\Theta^3$.} 
Our input is now the clean labeling $\Theta^2$.
Pick any $z_3\in\Gamma$ such that $\Theta_3^2(z_3)=q_3$.
We modify the labeling $\Theta^2$ on the ball
$B_{s_3}(G_{f(3)},z_3)$.
Let $\cal{A}^2=\{A^2_i\}_{i=1}^{\tau_2}$
be the set of all $\cF_2$-labeled $G_{f(2)}$-balls of radius $s_2$ in $\Theta^2$ up to
isomorphism. Let $\{w^2_i\}_{i=1}^{\tau_2}\subset B_{s_3}(G_{f(3)},z_3)\cap (\Theta_2^2)^{-1}(q_2)$ be a subset of 
vertices
such that:
\begin{itemize}
\item For any $1\leq i\leq \tau_2$,
$$\frac{1}{3} s_3\leq d_{G_{f(3)}}(z_3,w^2_i)\leq \frac{2}{3} s_3.$$
\item For any $1\leq i\neq j\leq \tau_2$,
$$d_{G_{f(3)}}(w^2_i,w^2_j)>20 r_2\,.$$
\end{itemize}
\noindent
Then for each $1\leq i \leq \tau_2$
pick $z^2_i\in (\Theta_2^2)^{-1}(q_2)$ such that
$B_{2r_2}^{\Theta^2,2}(G_{f(2)},z^2_i)$ fully contains an isomorphic copy of the ball 
$A^2_i$. 
Then using Proposition \ref{supersize}, for each $1\leq i \leq \tau_2$ we 
simultaneously patch the
balls $B_{3r_2}^{\Theta^2,2}(G_{f(2)},z^2_i)$ into $\Theta^2$ around $w^2_i$. 
Hence, we obtain a modified clean labeling $\hat{\Theta}^3:\Gamma\to\cF$.
So far, the construction of the $3$-codeball was identical to 
the one of $B_2$. Now, we need some further considerations.
It is possible that after the patching we will have some $g\in (\hat{\Theta}^3)^{-1}(q_2)$ such that
the ball $B_{s_2}^{\hat{\Theta}^3,2}(G_{f(2)},g)$ is not isomorphic to the codeball $B_2$. Using Proposition \ref{patching}, we patch the
ball $B^{\Theta^2,2}_{s_2+r_1}(G_{f(2)},z_2)$ into $\hatt^3$ around all such elements $g$ as above to obtain the clean labeling
$\hat{\Theta}^3_1$. We call the $\cF_3$-labeled $G_{f(3)}$-ball $B_3=B_{s_3}^{\hat{\Theta}_1^3,3}(G_{f(3)},z_3)$ the {\bf $3$-codeball}.
Now we construct the clean labeling $\ovet^3$ by patching the ball
$B^{\hatt_1^3,3}_{s_3+r_2}(G_{f(3)},z_3)$ into $\Theta^2$ around all vertices $y\in (\Theta^2_3)^{-1}(q_3)$.
Again, due to the patchings, it is possible that for some $g\in (\ovet^3)^{-1}(q_2)$,
$B_{s_2}^{\ovet^3,2}(G_{f(2)},g)$ is not isomorphic to the $2$-codeball $B_2$.
We obtain the clean labeling $\Theta^3$ by patching the ball $B_{s_2+r_1}^{\Theta^2,2}(G_{f(2)},z_2)$ 
into $\ovet^3$ around the vertices $g$ above.
Then
\begin{itemize}
\item For all $z\in(\Theta^3_3)^{-1}(q_3)$ we have
$B_{s_3}^{\Theta^3,3}(G_{f(3)},z)\cong B_3$.
\item For all $z\in(\Theta^3_2)^{-1}(q_2)$ we have
$B_{s_2}^{\Theta^3,2}(G_{f(2)},z)\cong B_2$.
\end{itemize}

\section{Minimal Cantor labelings III. (the induction)} \label{induction}
Let us suppose that we have already constructed the clean labelings \\ $\Theta^1,\Theta^2,\dots,\Theta^n$
satisfying the following properties.
\begin{itemize}
\item If $j\leq k \leq n$ and $(\Theta^k_j)(z)=q_j$, then the $\cF_j$-labeled $G_{f(j)}$-ball \\ $B_{s_j}^{\Theta^k,j}(G_{f(j)},z)$ is
isomorphic to a given $\cF_j$-labeled $G_{f(j)}$-ball $B_j$, the {\bf j-codeball}.
\item For $1< j\leq n$, let $\cal{A}^{j-1}=\{ A^{j-1}_i\}^{\tau_{j-1}}_{i=1}$ be the set of all $\cF_{j-1}$-labeled 
$G_{f(j-1)}$-balls of
radius $s_{j-1}$ in $\Theta^{j-1}$ up to isomorphism.  Then
for any $1\leq i \leq \tau_{j-1}$, the $j$-codeball $B_j$ fully contains a copy of $A_i^{j-1}$.
\end{itemize}
\noindent
Now we construct the clean labeling $\Theta^{n+1}$ in two rounds.

\noindent
{\bf Round one: The $n+1$-codeball}.   Pick $z_{n+1}\in\Gamma$ such that $\Theta^n_{n+1}(z_{n+1})=q_{n+1}$.
We modify the labeling $\Theta^n$ on the ball
$B_{s_{n+1}}(G_{f(n+1)}, z_{n+1})$.
Let $\cal{A}^n=\{A^n_i\}^{\tau_n}_{i=1}$ be
the set of all $\cF_n$-labeled $G_{f(n)}$-balls of radius $s_n$ in $\Theta^n$ up to isomorphism.
Let $\{w_i^n\}^{\tau_n}_{i=1}\subset B_{s_{n+1}}(G_{f(n+1)}, z_{n+1})\cap (\Theta^n_n)^{-1}(q_n)$ be a subset of vertices
such that
\begin{itemize}
\item For any $1\leq i\leq \tau_n$,
$$\frac{1}{3} s_{n+1}\leq d_{G_{f(n+1)}}(z_{n+1},w^n_i)\leq \frac{2}{3} s_{n+1}.$$
\item For any $1\leq i\neq j\leq \tau_n$,
$$d_{G_{f(n+1)}}(w^n_i,w^n_j)>20r_n\,.$$
\end{itemize}
\noindent
By Rule \ref{rule1}, such system indeed exists.
Now, for each $1\leq i \leq \tau_n$
pick $z^n_i\in (\Theta_n^n)^{-1}(q_n)$ such that
$B_{2r_n}^{\Theta^n,n}(G_{f(n)},z^n_i)$ fully contains an isomorphic copy of the ball 
$A^n_i$. 
Finally,  using Proposition \ref{supersize}, for each $1\leq i \leq \tau_n$ we 
simultaneously patch the
balls $B_{3r_n}^{\Theta^n,n}(G_{f(n)},z^n_i)$ into $\Theta^n$ around $w^n_i$. 
Hence, we obtain a modified clean labeling $\hat{\Theta}^{n+1}:\Gamma\to\cF$.
Now we need a definition to describe our procedure.
Let $\Theta':\Gamma\to \cF$ be a clean labeling and
let $2\leq i \leq n$. Let
$\mbox{Bad}_i(\Theta')$ be the set of all $y\in\Gamma$ such that $ \Theta'_i(y)=q_i$ and $B^{\Theta',i}_{s_i}(G_{f(i)},y)$
is not isomorphic to the $i$-codeball $B_i$.
The {\bf $i$-repair} of $\Theta'$, $\Theta'_{<i>}$ is constructed in the following way.
Using Proposition \ref{patching}, we simultaneously patch the
$\cF_i$-labeled $G_{f(i)}$-ball $B_{s_i+r_{i-1}}^{\Theta^i,i}(G_{f(i)},z_i)$ into $\Theta'$ around all elements $y\in\,\mbox{Bad}_i(\Theta')$.
Here, $z_i\in\Gamma$ is an element for which  $B_{s_i}^{\Theta^i,i}(G_{f(i)},z_i)$ is isomorphic to the $i$-codeball $B_i$.
Clearly, $\mbox{Bad}_i(\Theta'_{<i>})$ is empty. Note however, that if
$\mbox{Bad}_j(\Theta')$ is empty for some $j<i$, it is possible that $\mbox{Bad}_j(\Theta'_{<i>})$ is non-empty.
Now, we construct inductively the clean labelings $\{\hatt^{n+1,i}\}^{n+1}_{i=2}$.
Let $\hatt^{n+1,n+1}=\hatt^{n+1}$ and for $2\leq i \leq n$, let
$\hatt^{n+1,i}=\hatt^{n+1,i+1}_{<i>}$.
\begin{definition} Let $2\leq i \leq n$. Let $i>j_l>j_{l-1}>\dots> j_1\geq 2$
be integers and for $1\leq k \leq l$ let $g_k\in\Gamma$ such that
if $k<l$ then the set
$$B_{s_{j_k}+5 r_{j_k-1}}(G_{f(j_k)},g_k)\cap
B_{s_{j_{k+1}}+5 r_{j_{k+1}-1}}(G_{f(j_{k+1})},g_{k+1})$$
\noindent
is non-empty. Then, the set  $\cal{L}=\cup_{k=1}^l B_{s_{j_k}+5r_{j_k-1}}(G_{f(j_k)},g_k)$
is called an $i$-chain of length $l$.
\end{definition}
\begin{lemma} \label{chain}
The $G_{f(i)}$-diameter of an $i$-chain $\cal{L}$ is less than $\frac{1}{10} r_{i-1}$.
\end{lemma}
\proof
We have that
$$\mbox{diam}_{G_{f(i)}}(\cal{L})\leq \sum ^l_{k=1} \mbox{diam}_{G_{f(i)}}( B_{s_{j_k}+5r_{j_k-1}}(G_{f(j_k)},g_k))\,,$$
\noindent
hence the lemma follows from Rule \ref{rule2}. \qed
\begin{proposition}\label{1empty}
For any $2\leq i \leq n$,
$\mbox{Bad}_i(\hatt^{n+1,2})$ is empty.
\end{proposition}
\proof
Let $\hatt^{n+1,i}_i(w)=q_i$. It is enough to prove that
$B_{s_i}^{\hatt^{n+1,2},i}(G_{f(i)},w)\cong B_i$. 
Suppose that the ball $B_{s_i}^{\hatt^{n+1,2},i}(G_{f(i)}
,w)$ is not isomorphic
to the $i$-codeball $B_i$. Since
$B_{s_i}^{\hatt^{n+1,i},i}(G_{f(i)},w)\cong B_i$, there exists
$h\in\Gamma$ and $j<i$
such that
\begin{itemize}
\item $h\in \mbox{Bad}_{s_j}(\hatt^{n+1,j+1})$
\item $B_{s_j+5r_{j-1}}(G_{f(j)},h)\cap B_{s_i}(G_{f(i)},w)\neq 0$.
\end{itemize}
\noindent
Let $l$ be the largest integer such that we have an $i$-chain $\cal{L}$ of length $l$ satisfying
the following properties.
\begin{itemize}
\item $\cal{L} \cap B_{s_i}(G_{f(i)},w) \neq 0$.
\item $g_k\in\mbox{Bad}_{s_{j_k}}(\hatt^{n+1,j_k+1})$, where
$$\cal{L}=\cup_{k=1}^l B_{s_{j_k}+5r_{j_k-1}}(G_{f(j_k)},g_k)\,.$$
\noindent
So, by our previous observation we have that $l\geq 1$.
\end{itemize}
\begin{lemma}
The set
$$B_{s_{j_l}+5r_{j_l-1}}(G_{f(j_l)},g_l)\cap (\Gamma\backslash B_{s_i+\frac{1}{2}r_{i-1}}(G_{f(i)},w))$$ 
\noindent
is nonempty.
\end{lemma}
\proof
If $B_{s_{j_l}+5 r_{j_l-1}}(G_{f(j_l)},g_l)$ is contained in the ball $B_{s_i+\frac{1}{2}r_{i-1}}(G_{f(i)},w)$, then \\
$B_{s_{j_l}}^{\hatt^{n+1,i},j_l}(G_{f(j_l)},g_l)$ is isomorphic to the $j_l$-codeball $B_{j_l}$.
Thus, there is some $j_l<k<i$ and $g\in \Gamma$ such that $g\in \mbox{Bad}_k(\hatt^{n+1,k+1})$ and
the ball \\ $B_{s_k+5 r_{k-1}}(G_{f(k)},g)$ intersects the ball $B_{s_{j_l}+5 r_{j_l-1}}(G_{f(j_l)},g_l)$,
in contradiction with the maximality of $l$. Hence the lemma follows. \qed
\vskip 0.1in
\noindent
By Lemma \ref{chain}, $\mbox{diam}_{G_{f(i)}}(\cal{L})\leq \frac{1}{10} r_{i-1}$. Therefore, it is not possible
that both $\cal{L} \cap B_{s_i}(G_{f(i)},w) \neq \emptyset$ and $\cal{L}\cap (\Gamma\backslash B_{s_i+\frac{1}{2}r_{i-1}}(G_{f(i)},w))\neq \emptyset$ hold.
Therefore, our proposition follows. \qed
\vskip 0.1in
\noindent
We call the $\cF_{n+1}$-labeled $G_{f(n+1)}$-ball $B_{n+1}=B_{s_{n+1}}^{\hatt^{{n+1},2},n+1}(G_{f(n+1)},z_{n+1})$ 
the {\bf $n+1$-codeball}.
Now we construct the clean labeling
$\ovet^{n+1}$ by patching the ball $B_{s_{n+1}+r_n}^{\hatt^{{n+1},2},n+1}(G_{f(n+1)},z_{n+1})$ into $\Theta^n$ around
all the vertices $y\in (\Theta_{n+1}^n)^{-1}(q_{n+1})$.
\vskip 0.1in
\noindent
{\bf Round two: The construction of $\Theta^{n+1}$.} Again, it is possible that
for some $2\leq i \leq n$ the set
$\mbox{Bad}_i(\ovet^{n+1})$ is non-empty.
Now, we inductively construct the clean labelings $\{\ovet^{n+1,i}\}^{n+1}_{i=2}$.
Let $\ovet^{n+1,n+1}=\ovet^{n+1}$ and for $2\leq i \leq n$, let
$\ovet^{n+1,i}=\ovet^{n+1,i+1}_{<i>}$.
Copying the proof of Proposition \ref{1empty}, we can prove the following proposition.
\begin{proposition} \label{2empty}
For any $2\leq i \leq n$,
$\mbox{Bad}_i(\ovet^{n+1,2})$ is empty.
\end{proposition}
\noindent
Let $\Theta^{n+1}=\ovet^{n+1,2}$. Let us summarize what we have already proved about the clean labeling $\Theta^{n+1}$.
\begin{proposition}\label{summ}\mbox{} \begin{itemize}
\item For all $1< j \leq n+1$ and $z\in(\Theta^{n+1}_j)^{-1}(q_j)$ we have
$B_{s_j}^{\Theta^{n+1},j}(G_{f(j)},z)\cong B_j.$
\item For any $1< j\leq n+1$ and $1\leq i \leq \tau_{j-1}$, the
$j$-codeball $B_j$ fully contains a copy of the ball $A_i^{j-1}$.\end{itemize}
\end{proposition}
\noindent
Finally, we have the following lemma.
\begin{lemma} \label{tavol}
Let $z\in\Gamma$ such that
$\Theta^{n+1}(z)\neq \Theta^n(z)$. Then
$d_{G_{f(n+1)}}(z,e_\Gamma)\geq 5s_{n+1}$.
\end{lemma}
\proof 
By our construction if $\Theta^{n+1}(z)\neq \Theta^n(z)$ then there exists
some $y\in \Gamma$ such that 
\begin{itemize}
\item $\Theta^n_{n+1}(y)=q_{n+1}$ and there is an $l$-chain $\cal{L}$ for some $l\leq n$, 
intersecting $B_{2s_{n+1}}(G_{f(n+1)},y)$ for which $z\in B_{2s_{n+1}}(G_{f(n+1)},y)\cup \cal{L}$.
\end{itemize}
\noindent
Therefore, $z\in B_{5s_{n+1}}(G_{f(n+1)},y)$, hence the lemma follows from the definition of the clean labelings.
\qed
\section{Explicit construction of a free, minimal action}
Let $\{\Theta^n\}^\infty_{n=1}$ be the clean labelings constructed in the previous sections.
By Lemma \ref{tavol}, for any $y\in\Gamma$ there exists $n_y>0$ such that if $n\geq n_y$, then $\Theta^n(y)=\Theta^{n_y}(y)$. 
Hence, $\lim_{n\to\infty} \Theta^n(y)=\Theta^\infty(y)$ is a well-defined clean labeling.
\begin{lemma}
For any $i\geq 2$ and $z\in (\Theta^\infty_i)^{-1}(q_i)$, $B_{s_i}^{\Theta^\infty,i}(G_{f(i)},z)\cong B_i$.
\end{lemma}
\proof Let $n\geq i$ such that $\Theta^\infty(y)=\Theta^n(y)$ for
all $y\in B_{s_i}(G_{f(i)},z)$. 
Then, $(\Theta^n_i)(z)=q_i$. By Proposition \ref{summ}, $B_{s_i}^{\Theta^n,i}(G_{f(i)},z)\cong B_i$.
Hence our lemma follows. \qed
\begin{lemma} \label{synd}
Let $z\in\Gamma$, and $t,j\geq 1$ be integers. Then, there exists some $\rho>0$ such that
the $\cF_j$-labeled ball $B_t^{\Theta^\infty,j}(G_{f(j)},z)$ is $(\rho,j)$-syndetic in $\Theta^\infty$.
\end{lemma}
\proof
Let $q>j$ be an integer such that
the $\cF_j$-labeled ball $B_t^{\Theta^m,j}(G_{f(j)},z)$ is isomorphic to $B_t^{\Theta^\infty,j}(G_{f(j)},z)$, whenever $m\geq q$.
Let $n>q$ be an integer such that
$s_n>t$.
Then the $n+1$-codeball $B_{n+1}$ contains a copy of the $\cF_n$-labeled $G_{f(n)}$-ball $B_{s_n}^{\Theta^n,n}(G_{f(n)},z)$.
In particular, $B_{n+1}$ contains
a copy of the ball $B_t^{\Theta^n,j}(G_{f(j)},z)\cong B_t^{\Theta^\infty,j}(G_{f(j)},z)$.
The set $(\Theta^\infty_1)^{-1}(q_1)$ is a maximal $r_1$-sparse set, hence an $r_1$-net.
The set $(\Theta^\infty_2)^{-1}(q_2)$ is a $(\Theta^\infty_1)^{-1}(q_1)$-maximal $r_2$-sparse set, hence by Lemma \ref{net} an
$r_2+r_1$-net. Inductively, $(\Theta^\infty_{n+1})^{-1}(q_{n+1})$ is a $\sum_{i=1}^{n+1}r_i$-net, thus a $2r_{n+1}$-net.
So, for every $x\in\Gamma$ the ball
$B^{\Theta^\infty,n+1}_{3r_{n+1}}(G_{f(n+1)},x)$ contains a copy of
$B_t^{\Theta^\infty,j}(G_{f(j)},z)$.
That is, the $\cF_j$-labeled $G_{f(j)}$-ball $B_t^{\Theta^\infty,j}(G_{f(j)},z)$ is $(3r_{n+1},j)$-syndetic in $\Theta^\infty$. \qed.
\begin{theorem}\label{vegsoot}
The orbit closure of $\Theta^\infty$ in the Bernoulli shift space $\cF^\Gamma$ is free and minimal.
\end{theorem}
By Lemma \ref{synd}, the orbit closure of $\Theta^\infty$ is minimal. Since $\Theta^\infty$ is a clean labeling, it
is proper in the sense that for any $r>0$ there exists some $s_r>0$ such that
if $0<d_{G_r}(x,y)\leq r$, then
$$\Theta^{\infty}_{[s_r]}(x)\neq \Theta^{\infty}_{[s_r]}(y)\,.$$
\noindent
Therefore, by repeating the argument of the proof of Proposition \ref{free}, we immediately obtain that the action of $\Gamma$
on the orbit closure of $\Theta^\infty$ is free. \qed
\section{The dynamical version of clean labelings}
In this section we prove some propositions and lemmas that we need for the proof of Theorem \ref{fotetel2}.
For the whole section let the integers $\{s_i\}^\infty_{i=1}$, 
$\{r_i\}^\infty_{i=1}$ and the finite sets $\{F_m\}^\infty_{m=1}$ be as in the previous sections.
Our first proposition is about the existence of the dynamical analogue
of clean labelings.
\begin{proposition}\label{prelabeling}
Let $\alpha:\Gamma\actson \cC$ be a free continuous action as in the
previous sections. Then there exists a continuous map
$$\Sigma=\Sigma_{\cC}\times\prod_{m=1}^\infty \Sigma_m:\cC\to \cF$$
\noindent such that for any $x\in\cC$, $\Sigma_x(\gamma):=\Sigma(\alpha(\gamma)(x))$
defines an almost clean labeling of the group $\Gamma$ (see Section \ref{negy}).
\noindent
We call $\Sigma$ a {\bf dynamical clean labeling}.
\end{proposition}
\proof 
Before starting the proof let us fix some definitions and notations.
As in Section \ref{ketto}, we will consider the Borel graphs $\cG_r$ on the Cantor set $\cC$
associated to the free action $\alpha$.
If $U$ is a clopen subset of $\cC$, then we define a subset of our 
countable group $\Gamma$ by
$$O^U_x:=\{\gamma\in\Gamma\,\mid\,\alpha(\gamma)(x)\in U\}\,.$$
\noindent
Let $t>0$ be an integer. We call a clopen set $U\subset \cC$ {\bf $t$-separated} in the
Borel graph $\cG_r$ if
for any $x\neq y\in U$, $d_{\cG_r}(x,y)>t$. Clearly, $U$ is $t$-separated if and only if
for any $x\in\cC$ the set $O_x^U$ is $t$-sparse in the graph $G_r$.
\begin{lemma}\label{lemma81}
Let $U\subset \cC$ be a clopen set and $t>0$ be an integer.
Then there exists a clopen set $V\subset \cC$ such that
for any $x\in \cC$, the set $O^V_x$ is an $O^U_x$-maximal $t$-sparse set in the graph $G_r$.
\end{lemma}
\proof
By continunity and freeness, there exists some $s>0$
such that if for some $x,y\in \cC$ we have
$0<d_{\cG_r}(x,y)\leq t$, then $(x)_s\neq (y)_s$.
Let $a_1,a_2,\dots,a_{2^s}$ be an enumeration
of the set $\{0,1\}^{\{1,2,\dots,s\}}$.
For $1\leq i \leq 2^s$, let 
$$W_i:=\{x\in \cC\,\mid\,(x)_s=a_i\}\,.$$
Now let $V_1:=U\cap W_1$ and
$$Z_1=\cup_{d_{G_r}(\gamma,e_\Gamma)\leq t} \alpha(\gamma)V_1\,.$$
\noindent
Clearly, $Z_1$ is a clopen set. Let $Z_1^c$ denote the complement
of $Z_1$.
Now, let $V_2:=W_2\cap U \cap Z_1^c$. Then for any $x\in \cC$,
$O_x^{V_1\cup V_2}$ is a $t$-sparse set in $G_r$ contained in the set $O^U_x$.
Inductively, suppose that for $1\leq i \leq n< 2^s$,
\begin{itemize}
\item $V_i\subset W_i\cap U$ is a clopen set,
\item $Z_i=\cup_{d_{G_r}(\gamma,e_\Gamma)\leq t} \alpha(\gamma)(\cup^i_{j=1} V_j)\,,$
\item $V_{i+1}=W_{i+1}\cap U \cap Z^c_i$.
\end{itemize}
\noindent Then for any $1\leq i \leq n$, $O_x^{\cup^i_{j=1}V_j}$ is a $t$-sparse set in the graph $G_r$
which is contained in the set $O_x^U$.
We define
$$V_{n+1}:=W_{n+1}\cap U \cap Z_n^c\,.$$
\noindent
Then for any $x\in \cC$, $O_x^{\cup^{n+1}_{j=1}V_j}\subset O^U_x$ is a $t$-sparse set in the graph $G_r$.
We claim
that for any $x\in \cC$,
$O_x^{\cup^{2^s}_{j=1} V_j}$ is an $O_x^U$-maximal $t$-sparse set in $G_r$.
Clearly,
$O_x^{\cup^{2^s}_{j=1} V_j}\subset O_x^U$ and $O_x^{\cup^{2^s}_{j=1} V_j}$ is a $t$-sparse set in the graph 
$G_r$.
Suppose that $O_x^{\cup^{2^s}_{j=1} V_j}$ is not an $O_x^U$-maximal $t$-sparse set in the graph $G_r$.
Then there exists $\delta\in\Gamma_r$ such that
$\alpha(\delta)(x)\in U$ and 
\begin{equation} \label{mobi}
d_{G_r}(\delta, O_x^{\cup^{2^s}_{j=1} V_j})>t.
\end{equation}
Let $(\alpha(\delta)(x))_{s}=a_i$.
Then by \eqref{mobi},
$$\alpha(\delta)(x)\in W_i\cap U \cap Z_{i-1}^c=V_i$$
\noindent
in contradiction with our assumption. \qed
\begin{lemma} \label{lemma82}
Let $V$ be an $r_m$-separated clopen set in the Borel graph $\cG_{f(m)}$. Then there exists
a continuous function $\phi:\cC\to F_m$ such that 
\begin{itemize}
\item $V=\phi^{-1}(q_m)$.
\item For any $q\in F_m$, the
clopen set $\phi^{-1}(q)$ is $r_m$-separated in the Borel graph $\cG_{f(m)}$.
\end{itemize}
\end{lemma}
\proof
Let $V_1:=V$. Using Lemma \ref{lemma81}, we pick $V_2$
in such a way that for any $x\in \cC$, $O_x^{V_2}$ is an $O_x^{V_1^c}$-maximal $r_m$-sparse
set in the graph $G_{f(2)}$.
Inductively, suppose that the clopen sets $V_1,V_2,\dots, V_i$ have already been picked. Then, let $V_{i+1}\subset \cC$ be a clopen
set such that
for any $x\in \cC$, $O_x^{V_{i+1}}$ is an $O_x^{(\cup^i_{j=1} V_j)^c}$-maximal
$r_m$-sparse set in the graph $G_{f(m)}$.
It is enough to show that $\cup_{i=1}^{|F_m|} V_i=\cC$.
Suppose that $x\in \cC\backslash \cup_{i=1}^{|F_m|} V_i$. By Rule \ref{rule3}, there exists $1\leq i \leq |F_m|$, such that
for any $\gamma$ for which $d_{G_{f(m)}}(\gamma,e_\Gamma)\leq r_m$ holds, we have that  $\alpha(\gamma)(x)\notin V_i$.
However, this is in contradiction with the definition of the set $V_i$. \qed
\vskip 0.2in
\noindent
Now we finish the proof of our proposition. First, using Lemmas \ref{lemma81} and \ref{lemma82} we construct
a continuous map $\Sigma_1:\cC\to F_1$ such that for any $q\in F_1$ the set $(\Sigma_1)^{-1}(q)$ is $r_1$-separated in the Borel graph $\cG_{f(1)}$. Then,
using the lemmas we construct $\Sigma_2:\cC\to F_2$ such that
for any $q\in F_2$, $(\Sigma_2)^{-1}(q)$ is $r_2$-separated in $\cG_{f(2)}$ and
for any $x\in\cC$, $O_x^{(\Sigma_2)^{-1}(q_2)}$ is a $(\Sigma_1)^{-1}(q_1)$-maximal $r_2$-sparse set in the graph $G_{f(2)}$. Inductively, we construct $\Sigma_m$ for all $m\geq 1$. \qed
\vskip 0.2in
\noindent
{\bf Remark:} The inductive part in the proof of the previous proposition is the reason that in the definition of clean labelings  we required
$(\Theta_m)^{-1}(q_m)$ to be a $(\Theta_{m-1})^{-1}(q_{m-1})$-maximal $r_m$-sparse subset instead of
being a maximal $r_m$-sparse subset that is contained in $(\Theta_{m-1})^{-1}(q_{m-1})$.
\vskip 0.2in
\noindent
The following two lemmas are straightforward consequences of the definitions.
\begin{lemma} \label{lemma91}
Let 
$$\Sigma'=\Sigma'_{\cC} \times \prod^\infty_{m=1} \Sigma'_m:\cC\to \cF$$
and 
$$\Theta'=\Theta'_{\cC} \times \prod^\infty_{m=1} \Theta'_m:\cC\to \cF$$
be continuous functions.
Let $U\subset \cC$ be a $2t$-separated clopen set in the Borel graph $\cG_r$. 
Also, let $B^{\Theta,i}_t(G_r, e_\Gamma)$ be a $\cF_i$-labeled ball of radius $t$ in the graph $G_r$.
Let
$$\Sigma"=\Sigma"_{\cC} \times \prod^\infty_{m=1} \Sigma"_m:\cC\to \cC\times \prod^\infty_{m=1} F_m$$
be defined in the following way.
\begin{itemize}
\item If $d_{\cG_r}(x,U)>t$, then
$\Sigma"(x)=\Sigma'(x)$.
\item If there exists $y\in U$ and $\gamma\in\Gamma$, $d_{G_r}(\gamma,e_\Gamma)\leq t$ such that
$\alpha(\gamma)(y)=x$, then let
$\Sigma"_{[i]}(x)=\Theta'_{[i]}(\gamma)\,,(\Sigma"_{\cC}(x))_{i^+}=(\Sigma'_{\cC}(x))_{i^+}$ and for $j>i$, let $\Sigma"_j(x)=\Sigma'_j(x)$.
\end{itemize}
\noindent
Then $\Sigma"$ is a continuous map as well.
\end{lemma}
\begin{lemma} \label{lemma92}
Let $V\subset \cC$ be a $t$-separated clopen set. Let $s<t$ and suppose that
for some $\eps>0$ we have
$$|B_s(G_r,e_\Gamma)|\leq \eps |B_t(G_r,e_\Gamma)|\,.$$
\noindent
Let $U\subset \cC$ be the set of elements $x$ such that there exists $y\in V$ so that $d_{\cG_r}(x,y)\leq s$.
Then, for any Borel probability measure $\mu$ on $\cC$ that is invariant under the
action $\alpha$,
$$\mu(U)<\eps\,.$$
\end{lemma}
\section{The proof of Theorem \ref{fotetel2}} \label{end}
\noindent
Let $\zeta:\Gamma\actson \cC$ be a free continuous action.
Let $\{\Sigma^i:\cC\to \cF\}^\infty_{i=1}$ be a sequence of
dynamical clean labelings. We call $y\in \cC$ stable (with respect to the system $\{\Sigma^i\}^\infty_{i=1}$) if for
any element $x\in\cC$  in the $\Gamma$-orbit of $y$, there exists $n_x>0$ such that
if $n\geq n_x$ then $\Sigma^n(x)=\Sigma^{n_x}(x)$. So, for the $\Gamma$-invariant subset of stable points $y$,
$\lim_{n\to\infty}\Sigma^n(y)=\Sigma^\infty(y)$ exists.
\begin{proposition}\label{mainpropo}
There exists a sequence of dynamical clean labelings \\ $\{\Sigma^i\}^\infty_{i=1}$ such that
\begin{itemize}
\item For any Borel probability measure $\mu$ that is invariant under the action $\zeta$, the $\mu$-measure of
the stable points is $1$.
\item There exists a free and minimal Bernoulli subshift \\ $M\subset \cF^\Gamma$
such that for each stable point $y\in \cC$, 
$$\gamma \to \Sigma^\infty(\zeta(\gamma)(y))$$
\noindent
defines an injective equivariant Borel map from the Borel set of stable points to $M$.
\end{itemize}
\end{proposition}
\proof
We proceed very similarly as we did in the construction of $\Theta^\infty$.
The following lemma is an immediate corollary of Lemma \ref{lemma91}.
\begin{lemma} \label{lemma101}
Let $$\Sigma^a:\Sigma^a_\cC\times \prod^\infty_{m=1}\Sigma^a_m: \cC\to \cF$$ be a dynamical clean labeling and
$$\Theta^d: \Theta^d_\cC\times \prod^\infty_{m=1}\Theta^d_m:    \Gamma\to \cF$$ be an almost clean labeling.
Let $\Theta^d_n(x)=q_n$. Let $V\subset \cC$ be an $r_n$-separated clopen subset. For all $y\in V$
we patch the ball $B^{\Theta^d,n}_{s_n+r_{n-1}}(G_{f(n)},x)$ into $\Sigma^a_y$ around $y$. Then we obtain
a new dynamical clean labeling $\Sigma^b$.
\end{lemma}
\vskip 0.1in
\noindent
{\bf The construction of $\Sigma^2$.}
Let $\Sigma^1$ be a dynamical clean labeling such that for $x=(x_1,x_2,x_3,\dots)\in \cC$,
$\Sigma^1_{\cC}(x)=(y_1,y_2,y_3,\dots),$ where for $n\geq 0$, $y_{2^n}=x_1$ and if $2^n<a<2^{n+1}$, then $y_a=x_{a-2^n+1}$.
\noindent
Let $x\in\cC$ and $\Theta^1:=\Sigma^1_x$.
Let $\hat{\cal{A}}^1=\{\hat{A}^1_i\}_{i=1}^{\hat{\tau}_1}$
be the set of all $\cF_1$-labeled $G_{f(1)}$-balls of radius $s_1$ in the almost
clean labelings $\{\Sigma^1_w\}_{w\in\cC}$. Note
that $\hat{\cal{A}}^1$ might contain balls that are not in $\Theta^1$.
Let $z_2\in(\Theta^1_2)^{-1}(q_2)$. Proceeding in the same way as in Section \ref{codeball}
(using $\hat{\cal{A}}^1$ instead of $\cal{A}^1$), we construct the modified
labeling $\hatt^2:\Gamma\to \cF$. Now, our new $2$-codeball will be
$\hat{B}_2=B^{\hatt^2,2}_{s_2}(G_{f(2)},z_2)$. Then, we patch the ball
$B^{\hatt^2,2}_{s_2+r_1}(G_{f(2)},z_2)$ into $\Sigma^1$ around all vertices
$y\in(\Sigma^1)^{-1}(q_2)$ to obtain (Lemma \ref{lemma101}) the dynamical clean labeling $\Sigma^2$. 
\vskip 0,1in
\noindent
{\bf The induction.}
Let us suppose that we have already constructed the dynamical clean labelings \\ $\Sigma^1,\Sigma^2,\dots,\Sigma^n$
satisfying the following properties.
\begin{itemize}
\item If $j\leq k \leq n$, $x\in\cC$, $\Theta^k=\Sigma^k_x$ and $(\Theta^k_j)(z)=q_j$, then the $\cF_j$-labeled $G_{f(j)}$-ball $B_{s_j}^{\Theta^k,j}(G_{f(j)},z)$ is
isomorphic to a given $\cF_j$-labeled $G_{f(j)}$-ball $\hat{B}_j$, our new j-codeball.
\item For $1< j\leq n$, let $\hat{\cal{A}}^{j-1}=\{ \hat{A}^{j-1}_i\}^{\hat{\tau}_{j-1}}_{i=1}$ be the set of all $\cF_{j-1}$-labeled $G_{f(j-1)}$-balls of
radius $s_{j-1}$ in the almost clean labelings $\{\Sigma^{j-1}_w\}_{w\in\cC}$ up to isomorphism.  Then
for any $1\leq i \leq \hat{\tau}_{j-1}$, the $j$-codeball $\hat{B}_j$ fully contains a copy of $\hat{A}_i^{j-1}$.
\end{itemize}
\noindent
Again, we construct the dynamical clean labeling $\Sigma^{n+1}$ in two rounds.
\vskip 0.1in
\noindent
{\bf Round one: The new $n+1$-codeball.}
Let $x\in\cC$ and $\Theta^n=\Sigma^n_x$.
Also, let $z_{n+1}\in\Gamma$ such that
$\Theta^n_{n+1}(z_{n+1})=q_{n+1}$.
Let $\hat{\cal{A}}^n=\{\hat{A}^n_i\}^{\hat{\tau}_n}_{i=1}$ be
the set of all $\cF_n$-labeled $G_{f(n)}$-balls of radius $s_n$ in the almost clean labelings $\{\Sigma^n_w\}_{w\in\cC}$ up to isomorphism.
Now we construct the almost clean labeling $\hatt^{n+1,2}:\Gamma\to\cF$ exactly the same way as in Section \ref{induction}.
The $\cF_{n+1}$-labeled $G_{f(n+1)}$-ball $\hat{B}_{n+1}=B_{s_{n+1}}^{\hatt^{n+1,2},n+1}(G_{f(n+1)},z_{n+1})$ will be the new $n+1$-codeball. 
Finally, we construct the dynamical clean labeling $\overline{\Sigma}^{n+1}$ by patching the ball
$B_{s_{n+1}+r_n}^{\hatt^{n+1,2},n+1}(G_{f(n+1)},z_{n+1})$ into $\Sigma^n$ around all the vertices
$y\in (\Sigma^n_{n+1})^{-1}(q_{n+1})\,.$
\vskip 0.1in
\noindent
{\bf Round two: The construction of $\Sigma^{n+1}$.}
For a dynamical clean labeling $\Sigma'$, the set $\mbox{Bad}_i(\Sigma')$ is clopen so by Lemma \ref{lemma91} we can
construct the dynamical clean labelings $\overline{\Sigma}^{n+1,n+1}$, $\overline{\Sigma}^{n+1,i}=\overline{\Sigma}_{<i>}^{n+1,i+1}$ in the same
way as in Section \ref{induction}. Let $\Sigma^{n+1}=\overline{\Sigma}^{n+1,2}$ and let $Y\subset \cC$ be
the set of stable points with respect to the family of dynamical clean labelings $\{\Sigma^i\}^\infty_{i=1}$.
Then $\Sigma^\infty(y)=\lim_{n\to\infty}\Sigma^n(y)$ is a well-defined Borel map on the invariant Borel set $Y$.
By Theorem \ref{vegsoot}, for any $y\in Y$, the orbit closure of $\Theta^\infty_y$ in the Bernoulli shift
space $\cF^\Gamma$ is free and minimal. 
By our construction of the codeballs, if $y_1,y_2\in Y$, then for any $i\geq 1$, the set
of $\cF_i$-labeled $G_{f(i)}$-balls in $\Theta^\infty_{y_1}$ and $\Theta^\infty_{y_2}$ coincide. Hence the
orbit closures of $\Theta^\infty_{y_1}$ and $\Theta^\infty_{y_2}$ coincide as well. 
So, let $M$ be the orbit closure of $\Theta^\infty_y$ in $\cF^\Gamma$,
where $y\in Y$.
Notice that the map $\Phi:Y\to M$ defined by $\Phi(y)(\gamma)=\Sigma^\infty(\zeta(\gamma)(y))$ is injective.
Indeed, if $y_1\neq y_2\in Y$, then $\Sigma^1_\cC(y_1)$ and  $\Sigma^1_\cC(y_2)$ differ from each other in
infinitely many digits. Hence, $\Sigma^\infty(y_1)\neq \Sigma^\infty(y_2)$.
We finish the proof of our proposition by showing that for any Borel probability measure $\mu$ invariant under the action $\zeta$,
$\mu(Y)=1$.
Let $Q_n\subset \cC$ be the set of elements $w$ such that $\Sigma^{n+1}(w)\neq \Sigma^n(w)$.
By the Borel-Cantelli Lemma, it is enough to prove that
\begin{equation} \label{vegevege}
\sum^\infty_{n=1} \mu(Q_n)<\infty.
\end{equation}
Let $z\in Q_n$ and let
$V=(\Sigma^n_{n+1})^{-1}(q_{n+1})$. As we observed in the proof of Lemma \ref{tavol}, there
exists $y\in V$ and $\gamma\in B_{5s_{n+1}}(G_{f(n+1)},e_\Gamma)$ such that
$\zeta(\gamma)(y)=z$.
Note that $V$ is a $\frac{1}{10} r_{n+1}$-separated clopen set and by Rule \ref{rule2},
$$|B_{\frac{1}{10} r_{n+1}}(G_{f(n+1)}, e_\Gamma)|\geq 10^{n+1} |B_{5s_{n+1}}(G_{f(n+1)},e_\Gamma)|\,.$$
\noindent
Thus by Lemma \ref{lemma92} we have that $\mu(Q_n)<10^{-(n+1)}$, therefore \eqref{vegevege} holds.
This finishes the proof of our proposition. \qed.
\vskip 0.1in
\noindent
Now we can finish the proof of Theorem \ref{fotetel2}. Let $\zeta:\Gamma\actson \cC$ be the universal action of Theorem \ref{fotetel3}.
Let $\alpha:\Gamma\actson M$ be the Bernoulli subshift and $\{\Sigma^i\}^\infty_{i=1}$ be the sequence of dynamical
clean labelings in Proposition \ref{mainpropo}. Let $\beta:\Gamma\actson X$ be 
a Borel action,
let $\nu$ be a Borel probability measure on $X$ invariant under $\beta$, and finally let $\Psi_\beta:X\to \cC$ be 
the equivariant map in Theorem \ref{fotetel3}.
Then $(\Psi_\beta)_*(\nu)$ is a Borel probability measure invariant under the action $\zeta$. Let $Y\subset \cC$ be the set of stable points with respect
to the system $\{\Sigma^i\}^\infty_{i=1}$. Then, $\Phi(y)(\gamma)= \Sigma^\infty(\zeta(\gamma)(y))$ defines an injective, equivariant Borel map $\Phi:Y\to M$. 
Let $X'=(\Psi_\beta)^{-1}(Y).$
Then by Proposition \ref{mainpropo}, $\nu(X')=1$ and $\Phi\circ \Psi_\beta:X'\to M$ is an injective, equivariant Borel map. Hence, our theorem follows.
\qed

\vskip 0.2in
\noindent
Finally, we prove Corollary \ref{teteloxt}.
Let $\beta:\Gamma\actson \cC$ be the action in Theorem \ref{fotetel2}.
First let us suppose that $\Gamma$ is amenable.
Let $\alpha:\Gamma\actson (X,\mu)$ be an ergodic action and $\Phi_\alpha(X')\to \cC$ be an equivariant
injective map. Then the Kolmogorov-Sinai entropies of $\alpha$ and $(\Phi_\alpha)_*(\mu)$ coincide. Also, $(\Phi_\alpha)_*(\mu)$ is ergodic.
Since for any $0<c\leq \infty$ we have an ergodic action with Kolmogorov-Sinai entropy $c$, our corollary follows for amenable groups.
Now let us suppose that $\Gamma$ is non-amenable.
Again, we need to prove that there exist more than one non-isomorphic essentially free ergodic actions of $\Gamma$. 
However, by the famous theorem of Epstein
every non-amenable group has uncountably many pairwise orbit-inequivalent actions \cite{Eps} so our corollary follows for non-amenable groups as well.

\end{document}